\documentclass{amsproc}

\usepackage{amsfonts, amssymb, amscd}
\usepackage[all]{xy}

\copyrightinfo{2009}{American Mathematical Society}

\newcommand{\la}{\lambda}
\newcommand{\al}{\alpha}
\newcommand{\be}{\beta}

\newcommand{\f}{\varphi}

\newcommand{\abs}[1]{\vert #1\vert}

\newcommand{\Ll}{\operatorname{Lie}}

\newcommand{\CC}{\mathbb{C}}

\newcommand{\RR}{\mathbb{R}}
\newcommand{\ZZ}{\mathbb{Z}}
\newcommand{\NN}{\mathbb{N}}

\numberwithin{equation}{section}

\def\eqref#1{(\ref{#1})}

\newcommand{\arrow}{{\:\longrightarrow\:}}
\newcommand{\Z}{{\mathbb Z}}
\newcommand{\C}{{\mathbb C}}
\newcommand{\R}{{\mathbb R}}
\newcommand{\Q}{{\mathbb Q}}

\def\1{\sqrt{-1}\:}
\newcommand{\restrict}[1]{{\left|_{{\phantom{|}\!\!}_{#1}}\right.}}
\newcommand{\cntrct}                
{\hspace{2pt}\raisebox{1pt}{\text{$\lrcorner$}}\hspace{2pt}}

\renewcommand{\bar}{\overline}
\renewcommand{\phi}{\varphi}
\renewcommand{\epsilon}{\varepsilon}
\renewcommand{\geq}{\geqslant}
\renewcommand{\leq}{\leqslant}

\newcommand{\id}{\operatorname{\text{\sf id}}}

\newcommand{\Hom}{\operatorname{Hom}}

\newcommand{\rk}{\operatorname{rk}}

\renewcommand{\Im}{\operatorname{Im}}

\newtheorem{theorem}{Theorem}[section]
\newtheorem{proposition}{Proposition}[section]
\newtheorem{lemma}[theorem]{Lemma}
\newtheorem{corollary}{Corollary}[section]
\newtheorem{conjecture}{Conjecture}
\newtheorem{problem}{Open Problem}

\theoremstyle{definition}
\newtheorem{definition}[theorem]{Definition}

\theoremstyle{remark}
\newtheorem{remark}[theorem]{Remark}

\numberwithin{equation}{section}

\begin{document}

\title[Locally conformally K\"ahler manifolds]{A report on locally conformally K\"ahler manifolds}


\author{Liviu Ornea}

\address{University of Bucharest, Faculty of Mathematics, 14
Academiei str., 70109 Bucharest, Romania \emph{and} 
Institute of Mathematics ``Simion Stoilow" of the Romanian Academy, 
21, Calea Grivitei str.
010702-Bucharest, Romania }
\curraddr{}
\email{lornea@gta.math.unibuc.ro, liviu.ornea@imar.ro}
\thanks{L.O. is partially supported by a CNCSIS PNII IDEI Grant nr. 529/2009}

\author{Misha Verbitsky}
\address{Institute of Theoretical and
Experimental Physics \\
B. Cheremushkinskaya, 25, Moscow, 117259, Russia }
\curraddr{}
\email{verbit@verbit.ru}
\thanks{}

\subjclass[2000]{Primary 53C55}

\date{\today}

\dedicatory{To J. C. Wood at his sixtieth birthday}

\begin{abstract}
We present an overview of recent results in locally conformally K\"ahler geometry, with focus on the topological properties which obstruct the existence of such structures on compact manifolds. 
\end{abstract}

\maketitle

\section{Locally conformally K\"ahler manifolds}

Locally conformally K{\"a}hler (LCK) geometry is concerned with  
complex manifolds of complex dimension at least two 
admitting a K{\"a}hler covering
with deck transformations acting by holomorphic homotheties with respect to the K\"ahler metric. 

We shall
usually denote with $M$ the LCK manifold, with $(J,g)$ its Hermitian structure, with $\Gamma\rightarrow \tilde M\rightarrow M$ the K\"ahler covering and with $\tilde\omega$ the K\"ahler form on the covering.

Directly from the definition, one obtains the existence of  an associated character 
$$\chi:\Gamma\rightarrow \mathbb{R}^{>0}, \quad \displaystyle\chi(\gamma)=\frac{\gamma^*\tilde \omega}{\tilde\omega}.$$

This already puts some restrictions on $\pi_1(M)$. Others, more precise ones, will be obtained further.

Sometimes, the couple $(\Gamma, \tilde M)$ is called a {\bf presentation} of the LCK manifold $M$. Here, $\tilde M$ is understood as a  K\"ahler manifold tgether with a group of holomorphic homotheties (called a homothetic K\"ahler manifold). The idea is that, as on $M$ the metric can move in a conformal class, on the covering, the K\"ahler metric is not fixed but can be changed homothetically. Obviously, the same LCK manifold can admit many presentations and one can choose a minimal one and a maximal one (corresponding to the simply connected $\tilde M$). However, {\em the rank of the image of $\Gamma$ in $\mathbb{R}^{>0}$ is constant}; it will be denoted $\rk(M)$. Clearly, $\rk(M)\leq b_1(M)$ (see \cite{gopp}).

An equivalent definition  - historically, the first one -, at the level of the manifold itself, requires 
the existence of an open covering $\{U_\al\}$ with  local K\"ahler
metrics $g_\al$ subject to the condition that on overlaps $U_\al\cap U_\be$, these
local K\"ahler metrics are homothetic: $g_\al=c_{\al\be}g_\be$. The cocycle $\{c_{UV}\}$ is represented by a closed one-form $\theta$. 
Locally, $\theta\restrict{U_\al}=df_\al$ and the metrics
$e^{f_\al}g_\al$ glue to a global metric whose associated two-form
$\omega$ satisfies the integrability condition
$d\omega=\theta\wedge\omega$, thus being locally conformal with the
K\"ahler metrics $g_\al$. Here $\theta$ is a closed 1-form
on $M$, called {\bf the Lee form}. This gives
another definition of an LCK structure (motivating also the name), which
will be used in this paper.

\begin{definition}
Let $(M, \omega)$ be a Hermitian manifold,
$\dim_\C M >1$, with $d\omega = \theta\wedge \omega$,
where $\theta$ is a closed 1-form. Then
$M$ is called {\bf a locally conformally K{\"a}hler} (LCK)
manifold.
\end{definition}

\begin{remark}\label{rem1}
$i)$ Some authors include the K\"ahler manifolds  as particular LCK manifolds. Although this is a legitimate choice, we prefer the dichotomy LCK {\em versus} K\"ahler, and hence we always assume that the LCK manifolds we work with are of non-K\"ahler type. Due to a result in \cite{Va_80}, namely: {\em A compact locally conformally K\"ahler manifold which admits some K\"ahler
metric, or, more generally, which satisfies the $\partial\bar\partial$-Lemma, is globally conformally K\"ahler}, it is enough to assume $[\theta]\neq 0\in H^1(M,\R)$.

 $ii)$ The equation $d\omega=\theta\wedge\omega$ makes sense also in absence of a complex structure, leading to the notion of {\bf locally conformally symplectic manifold} (LCS). There is a great number of papers on this topic, among the  authors of which we cite:  A. Banyaga, G. Bande, S. Haller,  D. Kotschick, A. Lichnerowicz, J.C. Marrero, I. Vaisman etc. Hence, any LCK structure underlies a LCS structure. Nothing is known on the converse. The corresponding question regarding the relation symplectic {\em versus} K\"ahler was since long solved by Thurston, \cite{th}. We still do not know if any (compact) LCS manifold admits an integrable, compatible complex structure which makes it a LCK manifold or not. The difficulty might come from the fact that the topology of a LCK manifold is not controlled. We believe the answer should be negative and hence we propose:
\end{remark}

\begin{problem}
 Construct a compact LCS manifold which admits no LCK metric.
\end{problem}

The Lee form, which is the torsion of the Chern connection (see \cite{paul}), can also be interpreted in terms of presentations as follows. Abelianize the Serre sequence of $\Gamma\rightarrow \tilde M\rightarrow M$ to get:
$$H_1(K,\ZZ)\rightarrow H_1(M,\ZZ)\rightarrow \Gamma^{ab}\rightarrow 1.$$ 
Then apply $\Hom(\cdot, \ZZ)$ to obtain:
$$0\rightarrow \Hom(\Gamma^{ab}, \ZZ)\rightarrow H^1(M,\ZZ)\rightarrow H^1(K,\ZZ).$$
Tensoring with $\otimes_\ZZ\RR$, exactness is conserved (as $\RR_{\slash\ZZ}$ is flat) and one arrives at:
$$0\rightarrow \Hom_\ZZ(\Gamma^{ab},\RR)\stackrel{i}{\rightarrow} H^1_{DR}(M)\rightarrow H^1_{DR}(K).$$
Then  $i(\chi)=[\theta]$, as proven in \cite{pv}.

\smallskip

We refer to \cite{drag} and \cite{lsurv} for an overview of this geometry. Here we focus on our recent results and on related ones.

\smallskip

The following notion, coming from conformal geometry, is crucial for the way we understand LCK geometry:

\begin{definition}
Let  $(M,\omega, \theta)$ be an LCK manifold,
and $L$ the trivial line bundle, associated
to the representation $\mathrm{GL}(2n,\mathbb{R})\ni A \mapsto | \det
A|^{\frac{1}{n}}$, with flat connection
defined as $D:= d + \theta$. Then $L$ is called {\bf the
weight bundle} of $M$. 

Its holonomy coincides with the character $\chi:
\pi_1(M) \arrow \R^{>0}$ whose
image is called {\bf the monodromy group of $M$}.
\end{definition}

We shall denote with the same letter,
$D$, the corresponding  covariant derivative on $M$. It is precisely the Weyl covariant derivative associated to $\nabla=\nabla^g$ and $\theta$, uniquely defined by the conditions:
$$D J=0,\quad\quad 
D g=\theta\otimes g.$$ 

The complexified weight bundle will also be denoted $L$. As $d\theta=0$, $L$ is flat and defines a local system and hence one can compute its cohomology. 

On the other hand, in LCK geometry, one tries to work on the K\"ahler covering. But there, the interesting tensorial objects, in particular differential forms $\al$, are the ones satisfying:  $\gamma^*\al=\chi(\gamma)\al$ for every $\gamma\in \Gamma$. We call such forms {\bf automorphic}. 

The advantage of using the weight bundle is that automorphic objects on $\tilde M$ are regarded as objects on $M$ with values in $L$.

\subsection{Examples}

\subsubsection{Diagonal Hopf manifolds.} (\cite{go}, \cite{_Kamishima_Ornea_}, \cite{ve}.) Let $H_A:=(\C^n\setminus \{0\})/ \langle A\rangle$ with
$A=\text{diag}(\al_i)$ endowed with: 
\begin{itemize}
\item Complex structure as the 
projection of the standard one of $\C^n$.
\item LCK metric constructed as follows: 

\noindent Let $C>1$ be a constant and
$$\f(z_1,\ldots,z_n)=\sum \abs{z_i}^{\be_i}, \quad
\be_i=\log_{\abs{\al_i}^{-1}}C$$
a potential on $\C^n$.

Then  $A^*\f=C^{-1}\f$ and hence: $\displaystyle \Omega = \sqrt{-1} \partial\bar\partial \f$ is K\"ahler and
$\displaystyle \Gamma\cong \ZZ$ acts by holomorphic homotheties with
respect to it.
\end{itemize}
Note that the Lee field: $\displaystyle \theta^\sharp=-\sum z_i\log
\abs{\al_i}\partial z_i$ is parallel.

It is also important to observe that the LCK metric here is constructed out of an automorphic potential. The construction will be extended to non-diagonal Hopf manifolds.

\smallskip

\subsubsection{Compact complex surfaces.} Belgun, \cite{be}, gave the complete list  of compact complex surfaces which admit metrics with parallel Lee form ($\nabla\theta=0$), being, in particular, LCK. Such metrics are called {\bf Vaisman} and will be treated separately, in section \ref{vai} (see Theorems \ref{belg},  \ref{elsurf}).

Recently, Fujiki and Pontecorvo constructed LCK metrics on parabolic and hyperbolic Inoue surfaces. These examples are also bihermitian and hence related to generalized K\"ahler geometry. We also note that in \cite{ad}, the LCK metric of the diagonal Hopf surface $g_{GO}$ found \cite{go} was deformed to a family of bihermitian
metrics $(g_t, J, J^t_−)$ with $J^t_−
= \f^∗_t (J)$, where $\f_t$ is a path of diffeomorphisms;
as $t\rightarrow 0$, $J^t\rightarrow J$ and $g_t/t \rightarrow g_{GO}$.

More generally, Brunella, \cite{br2}, proved that all surfaces with global spherical shells, also known as Kato surfaces (as the previous mentioned parabolic and hyperbolic Inoue surfaces are) do admit LCK metrics. Previously he constructed families of LCK metrics only on Enoki surfaces, \cite{br1}.

On the other hand, Belgun also proved in \cite{be} that a certain type of Inoue surfaces does not admit any LCK metric. As these surfaces are deformations of other Inoue surfaces with LCK metric (found in \cite{_Tricerri_}), this proves, in particular, that, unlike the K\"ahler class, {\em the LCK class is not stable at small deformations}. By contrast, the LCK class share with the K\"ahler one the stability to blowing up points, \cite{_Tricerri_}, \cite{vuli}.

\smallskip

\subsubsection{Oeljeklaus-Toma manifolds,} \cite{ot}. Let $K$ be an algebraic number field of degree $n:=(K:\Q)$. Let then $\sigma_1,\ldots, \sigma_s$ (resp. $\sigma_{s+1},\ldots,\sigma_n$) be the real (resp. complex) embeddings of $K$ into $\CC$, with $\sigma_{s+i}=\bar\sigma_{s+i+t}$, for $1\leq i\leq t$. Let  $\mathcal{O}_K$ be the ring of algebraic integers of $K$. Note that for any $s,t \in \NN$, there exist algebraic number fields with precisely $s$ real and $2t$ complex embeddings.

Using the embeddings $\sigma_i$, $K$ can be embedded in $\C^m$, $m=s+t$, by 
$$\sigma:K\rightarrow \C^m, \quad \sigma(a)=(\sigma_1(a),\ldots,
\sigma_m(a)).$$
This embedding extends to $\mathcal{O}_K$ and $\sigma(\mathcal{O}_K)$ is a lattice of rank $n$ in $\C^m$. This gives rise to a 
properly discontinuous action of $\mathcal{O}_K$ on $\C^m$. On the other hand, $K$ itself acts on $\C^m$ by
$$(a,z)\mapsto (\sigma_1(a)z_1,\ldots, \sigma_m(a)z_m).$$
Note that if $a\in \mathcal{O}_K$, $a\sigma(\mathcal{O}_K)\subseteq \sigma(\mathcal{O}_K)$. Let now $\mathcal{O}_K^*$ be the group of units in $\mathcal{O}_K$ and set
$$\mathcal{O}_K^{*,+}=\{a\in \mathcal{O}_K^* \mid \sigma_i(a)>0, \, 1\leq i\leq s\}.$$
The only torsion elements in the ring $\mathcal{O}_K^*$ are $\pm 1$, hence the Dirichlet units theorem asserts the existence of a free Abelian group $G$ of rank $m-1$ such that $\mathcal{O}_K^* =G\cup (-G)$. Choose $G$ in such a  away that it contains $\mathcal{O}_K^{*,+}$ (with finite index). Now  $\mathcal{O}_K^{*,+}$acts multiplicatively on $\C^m$ and, taking into account also the above additive   action, one obtains a free action of the semi-direct product $\mathcal{O}_K^{*,+}\ltimes \mathcal{O}_K^*$ on $\C^m$ which leaves invariant $H^s\times \C^t$ (as above, $H$ is the open upper half-plane in $\C$). The authors then show that it is possible to choose a subgroup $U$ of $\mathcal{O}_K^{*,+}$ such that the action of $U\ltimes \mathcal{O}_K$ on $H^s\times \C^t$ be properly discontinuous and co-compact. Such a subgroup $U$ is called \emph{admissible} for $K$. The quotient 
$$X(K,U):=(H^s\times \C^t)/(U\ltimes \mathcal{O}_K)$$
is then shown to be a $m$-dimensional compact complex (affine) manifold, differentiably a fiber bundle over $(S^1)^s$ with fiber $(S^1)^n$.

For $t=1$, $X(K,U)$ admits LCK metrics.

Indeed, 
$$\f:H^s\times \C\rightarrow \RR, \quad \f=\prod_{j=1}^s\frac{i}{z_j-\bar z_j}+\abs{z_m}^2$$
is a K\"ahler potential on whose associated 2-form $i\partial\bar\partial \f$ the deck group acts by linear holomorphic homotheties. On the other hand, one sees that the potential itself is not automorphic (in particular, these manifolds cannot be Vaisman, see \S \ref{vai}).

A particular class of manifolds $X(K,U)$ is that of \emph{simple type}, when $U$ is not contained in $\ZZ$ and its action on  $\mathcal{O}_K$ does not admit a proper non-trivial invariant submodule of lower rank (which, as the authors show, is equivalent to the assumption that there is no proper intermediate field extension $\Q\subset K'\subset K$ with $U\subset \mathcal{O}_{K'}$). If $X(K,U)$ is of simple type, then 
$b_1(X(K,U))=s$ (a more direct proof than the original one appears in \cite{pv}), $b_2(X(K,U))=\genfrac{(}{)}{0pt}{0}{s}{2}$.
Moreover, the tangent bundle $T{X(K,U)}$ is flat and $\dim H^1(X(K,U), \mathcal{O}_{X(K,U)})\geq s$. In particular,  $X(K,U)$ are non-K\"ahler.

Observe that for $s=t=1$ and $U=\mathcal{O}_K^{*,+}$, $X(K,U)$ reduces to an Inoue surface $S_M$ with the metric given in \cite{_Tricerri_}.

Now, for $s=2$ and $t=1$, the six-dimensional $X(K,U)$ is of simple type, hence has the following Betti numbers: $b_0=b_6=1$, $b_1=b_5=2$, $b_2=b_4=1$, $b_3=0$. This disproves Vaisman's conjecture claiming that a compact LCK, non-K\"ahler,  manifold must have an odd odd Betti number.

These manifolds can be used to obtain examples of LCK structures with arbitrary rank (recall that $\rk (M)$ is the rank of $\chi(\Gamma)$ in $\R^{>0}$). Specifically:

\begin{theorem}\cite{pv}
 Let the number field $K$ admit exactly two non-real embeddings and $M=X(K,U)$. Then:

$i)$ If $n$ is odd (hence if $\dim_\C(M)$ is even), then $\rk (M)=b_1(M)$ ({\em i.e.} the rank is maximal).

$ii)$ If $n$ is even, then either $\rk M=b_1(M)$  or $\rk(M)=\displaystyle\frac{b_1(M)}{2}$; this last situation occurs if and only if $K$ is a
quadratic extension of a totally real number field.
\end{theorem}
Concrete examples of number fields which lead to $ii)$ above are also constructed in \cite{pv}.

\section{Locally conformally  K\"ahler manifolds with potential}

\begin{definition}\cite{_OV_Top_Potential_}
 $(M,J,g)$ is a LCK manifold with (automorphic) potential if $M$ admits a K\"ahler covering with automorphic potential.
\end{definition}

\begin{remark}
 The definition we gave in \cite{_OV:Potential_} was slightly more restrictive: we asked the potential to be a proper function ({\em i.e.} to have compact levels).  The properness of the potential is equivalent to the weight bundle having monodromy $\Z$. Later on, we proved in \cite{_ov:MN_}
that on any compact LCK manifold with automorphic potential, there exists another LCK metric 
with automorphic potential and monodromy $\Z$. The proof amounts to a deformation of the weight bundle together with its connection form.

However, we have strong reasons to believe that the deformation is not necessary:
\end{remark}

\begin{conjecture}
 Any compact LCK manifold with automorphic potential has monodromy $\Z$.
\end{conjecture}

The existence of a potential for the K\"ahler metric of the covering can be shown to be equivalent with the equation $(\nabla\theta)^{1,1}=0$, introduced in \cite{kok} under the name of pluricanonical K\"ahler-Weyl and studied also in \cite{_Kokarev_Kotschick:Fibrations_}.
 
\begin{proposition}\cite{_OV_Top_Potential_}\label{equiv_pot}
 $M$ admits a K\"ahler covering with automorphic potential  if and only if $(\nabla\theta)^{1,1}=0$.
\end{proposition}

For the proof, one first proves by direct computation that $(\nabla\theta)^{1,1}=0$ is equivalent with the equation:
$$d(J\theta) = \omega - \theta\wedge J\theta.$$
This can also be put in terms of Weyl connection as:
$$(D\theta)^{1,1} = (\theta\otimes \theta)^{1,1}-\frac 12 g.$$
Now, let $\tilde M$ be a covering of $M$ on which the pull-back
of $\theta$ is exact. Denote, for convenience, with the
same letters the pull-backs to $\tilde M$ of $\theta$,
$\omega$ and $D$. As locally $D$ is the Levi-Civita connection of the local K\"ahler metrics, its pull-back on $\tilde M$ is the Levi-Civita 
connection of the 
K\"ahler metric on $\tilde M$ globally conformal with
$\omega$. Then 
let $\psi:= e^{-\nu}$, where $d\nu = \theta$.
We have
\[  d d^c  \psi=
   -e^{-\nu} d d^c \nu + e^{-\nu} d\nu \wedge d^c \nu
   =  e^{-\nu} (d^c\theta + \theta\wedge J\theta)= \psi \omega,
\]
and hence the pluricanonical condition implies that $\psi$
is an automorphic potential for the K\"ahler metric $\psi\omega$. The
converse is true
by a similar argument.

\medskip

A second characterization can be given in terms of Bott-Chern cohomology. 
Let $\Lambda^{1,1}_{\chi,
  d}(\tilde M)$ be the space of closed, automorphic
$(1,1)$-forms on $\tilde M$, and $C^\infty_\chi(\tilde M)$
the space of automorphic functions on $\tilde M$. Then
\[
H^{1,1}_{BC}(M, L):= \frac{\Lambda^{1,1}_{\chi,
  d}(\tilde M)}{dd^c(C^\infty_\chi(\tilde M))}
\]
is the {\bf Bott-Chern group} of the LCK manifold (it is finite-dimensional and does not depend on the choice of the presentation). It is now clear that 

\begin{lemma}\cite{_ov:MN_}\label{potbc}
 {$M$ is LCK with potential if $[\Omega]=0\in H^{1,1}_{BC}(M, L)$.}
\end{lemma}

\smallskip

The main properties of LCK manifolds with automorphic potential are listed in the following:

\begin{theorem}\cite{_OV:Potential_}\label{emb}
$i)$ The class of compact LCK manifolds with
potential is stable to small deformations. 

$ii)$ Compact LCK manifolds with potential, of complex dimension at least $3$,
 can be holomorphically embedded in a (non-diagonal), 
Hopf manifold. 
\end{theorem}

{}From $i)$, it follows that the Hopf manifold $(\CC^N\setminus
0)/\Gamma$, with $\Gamma$ cyclic, generated by a
{\em non-diagonal} linear operator, is LCK with potential. 
This is the appropriate generalization of the (non--Vaisman)  non-diagonal  Hopf surface. Then $ii)$ says that the Hopf manifold plays in LCK geometry the r\^ole of the projective space in K\"ahler geometry.

\subsection{Vaisman manifolds}\label{vai}

Among the LCK manifolds with potential, a most interesting class is the Vaisman one. A Vaisman metric is a Hermitian metric with parallel Lee form. It can be easily seen that the K\"ahler metric of the covering has global automorphic potential $\f=\tilde\omega(\pi^*\theta,\pi^*\theta)$. 

The Lee field of a Vaisman manifold is Killing and, being parallel, it has constant length. Conversely, a LCK metric with Killing Lee field of constant length is Vaisman (see, {\em e.g.} \cite[Proposition 4.2]{drag}). On the other hand, it was proven in \cite[Proposition 6.5]{ve} that a complex compact submanifold of a compact Vaisman manifold must be tangent to the Lee field. In particular, the submanifold enherits a LCK metric whose Lee field is again Killing and of constant length. Hence:

\begin{proposition}
 Complex compact submanifolds of a compact Vaisman manifold are again Vaisman.
\end{proposition}

As the LCK metric of the diagonal Hopf manifold is Vaisman, this provides a wide class of examples.  On the other hand, on surfaces there exists the complete list of compact examples, see above.

On Vaisman manifolds, the vector field $\theta^\sharp$ is holomorphic and Killing, and hence it generates a totally geodesic, Riemannian, holomorphic foliation $\mathcal{F}$. When this is quasi-regular, one may consider the leaf space and obtain a fibration in elliptic curves over a K\"ahler orbifold. Similarly, when $\theta^\sharp$ has compact orbits, the leaf space is a Sasakian orbifold, \cite{bl}, over which $M$ fibers in circles. The two principal fibrations are connected by the Boothby-Wang fibration in a commutative diagram whose model is the classical Hopf fibering: 

$$
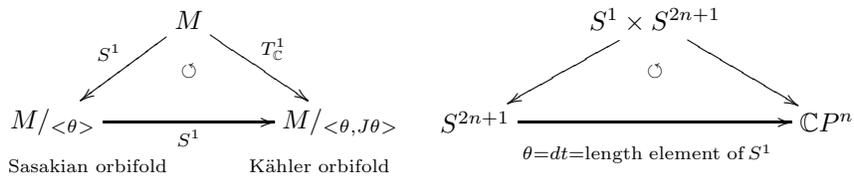
\underset{{\text{\!\!\!\scriptsize Sasakian orbifold\hspace{1cm} K\"ahler orbifold}}}
{
\xymatrix{&\ar[dl]_-{S^1}  M \ar@{}[d] |{\circlearrowleft} \ar[dr]^-{T^1_\mathbb{C}}&  \\
{M/_{<\theta>}}\ar[rr]_-{S^1} &&  M/_{<\theta, J\theta>}
}
}
\quad {\underset{{\theta=dt= \text{length element of}\, S^1}}{
\xymatrix{&\ar[dl]  S^1\times S^{2n+1} \ar@{}[d] |{\circlearrowleft} \ar[dr]&  \\
S^{2n+1}\ar[rr]  &&   \mathbb{C}P^n}
}
}
$$

\medskip

This is, in fact, the generic situation, because we proved in \cite{ov_ma1} that {\em the Vaisman structure of a compact manifold can always be deformed to a quasi-regular one}. 

\smallskip

{}From the above, it is clear that Vaisman structures may exist on the total space of some elliptic fibrations on compact K\"ahler manifolds. The precise statement is:

\begin{theorem}\cite{vu2}\label{vuli}
 Let $X,B$ be compact complex manifolds, $X\rightarrow B$ an elliptic
principal bundle with fiber $E$. If the Chern classes of this bundle  
are linearly independent in $H^2(B,\R)$, then $X$ carries no locally conformally
K\"ahler structure.
\end{theorem}
This contrasts with the case of an induced Hopf fibration over a projective manifold $B$, when one of the Chern classes vanishes.

\smallskip

For surfaces, we have a complete list of those who admit Vaisman metrics: 

\begin{theorem}\cite{be}\label{belg}
  Let $M$ be a compact complex surface with odd
$b_1$. Then $M$ admits a Vaisman metric if and only if
$M$ is an elliptic surface (a properly elliptic surface, a - primary or
secondary
- Kodaira surface, or an elliptic Hopf surface) or a diagonal Hopf surface.

\end{theorem}
 
Using the ``if'' part of this result we can prove:

\begin{theorem}\label{elsurf}
 Let $М$ be a minimal, 
non-K\"ahler compact surface,
which is not of class VII.
Then $M$ is a Vaisman elliptic
surface.
\end{theorem}

Indeed, recall that a compact complex surface 
surface is called {\bf class VII} if it 
has Kodaira dimension $-\infty$ and
$b_1(M)=1$. It is called {\bf minimal}
if it has no rational curves
with self-intersection $-1$. Now, from Kodaira's classification of surfaces,
it follows that the algebraic dimension
of $M$ is 1 (see {\em e.g.}  \cite[Theorem 5]{_Toma_}). 
Also from Kodaira's classification it follows that
$M$ is elliptic  \cite[Theorem 3]{_Toma_}. 
On the other hand, a non-K\"ahler compact complex surface has odd $b_1$ (\cite{_Buchdahl} and \cite{_Lamari}) 
It only remains to apply Belgun's result.

\smallskip 

The transversal K\"ahlerian foliation $\mathcal{F}$ permits the use of transversal foliations techniques (basic operators etc.) The following result concerning unicity of Vaisman structures was obtained this way:

\begin{theorem}\cite{_OV:EW_}
 Let $(M, J)$ be a compact complex manifold admitting a Vaisman
structure, and $V \in \Lambda^{n,n}(M)$ a nowhere degenerate, positive volume form. Then
$M$ admits at most one Vaisman structure with the same Lee class, such that the
volume form of the corresponding Gauduchon metric is equal to $V$.
\end{theorem}

Another recent application of this technique is the following:

\begin{theorem}\cite{op}
 Let $(M^{2m},g,J)$ be a compact Vaisman manifold. The metric $g$ is geometrically formal ({\em i.e.} the product of every harmonic forms is again harmonic) if and only if $b_p(M)=0$ for $2\leq p\leq 2m-2$ and $b_1(M)=b_{2m-1}(M)=1$, hence $M$ has the real homology of a Hopf manifold.
\end{theorem}

The connection between Vaisman and Sasakian geometries is clearly seen in:

\begin{theorem}\cite{ov_str}\label{ov_str}
 Compact Vaisman  manifolds are mapping tori over $S^1$ with Sasakian fibre. More precisely: 
 the universal cover $\tilde M$ is a metric cone $N \times \mathbb{R}^{>0}$, with $N$ compact Sasakian manifold
and the deck group is isomorphic with $\mathbb{Z}$, generated by
$(x, t)\mapsto (\lambda(x), t+q)$ for some $\lambda\in \mathrm{Aut}(N)$, $q\in \mathbb{R}^{>0}$.
\end{theorem}
 
This result was recently used to prove the following:

\begin{theorem}\cite{mo}\label{conf_hol}
 On compact Vaisman manifolds whose Weyl connection does not have holonomy in $\mathrm{Sp}(n)$ and which are not diagonal Hopf manifolds, $\mathfrak{conf}(M,[g])=\mathfrak{aut}(M)$.
\end{theorem}

Indeed, the statement follows from the fact that Killing fields with respect to the Gauduchon metric (and a Vaisman metric is Gauduchon) are holomorphic, \cite{mo}, and from the more general, referring to Riemannian cones:

\begin{theorem}\cite{mo}
 Let $(M,g):=(W,h)\times\RR/_{\{(x,t)\sim(\psi(x),t+1)\}}$, with $\psi\in \mathrm{Iso}(W,h)$, $W$ compact.
Then conformal vector fields on $(M,g)$ are
Killing.
\end{theorem}

For Vaisman manifolds, the conclusion of $ii)$ in Theorem \ref{emb} can be sharpened:

\begin{theorem}\cite{_OV:Potential_}\label{emb_vai}
 A compact complex manifold of dimension of least $3$ admits a Vaisman metric if and only if it admits a holomorphic embedding into a diagonal Hopf manifold.
\end{theorem}

Taking into account also the relation between Sasaki and Vaisman geometries, a first application of this Kodaira-Nakano type theorem was a corresponding embedding result in Sasakian geometry:

\begin{theorem}\cite{OV:emb_sas}
 A compact Sasakian manifold $M$ admits a CR-embed\-ding into a Sasakian manifold diffeomorphic to a sphere, and this embed\-ding is compatible with the respective Reeb fields. 
\end{theorem}

Moreover, we showed that this is the best result one may hope: assuming the existence of a model manifold in Sasakian geometry, analogue of the projective space in complex geometry, leads to a contradiction. A key point in the proof of the theorem was showing  that if $Z$ is a closed complex submanifold of a compact K\"ahler manifold $(M,\omega)$, $[\omega]\in H^2(M)$ is the K\"ahler class of $M$, and $\omega_0$ is a K\"ahler  form on $Z$ such that its K\"ahler class coincides with the restriction $[\omega]|_{Z}$, then there exists a K\"ahler form $\overline\omega\in[\omega]$ on $M$ such that $\overline\omega|_{Z}=\omega_0$. Recently, using a same type of argument, van Coevering gave a more direct proof of the embedding in \cite{cov}.

\smallskip

We also used Theorem \ref{emb_vai} to prove that, diffeomorphically, LCK with automorphic potential and Vaisman manifolds are the same:

\begin{theorem}\cite{_OV_Top_Potential_}
 Let $(M,\omega,\theta)$ be an LCK manifold
with potential with complex dimension at least $3$. Then there exists a  deformation of $M$ which
admits a Vaisman metric.
\end{theorem}
 For the proof, one considers a holomorphic embedding of $M$ in a Hopf manifold $H= (\C^N\setminus \{0\})\langle A\rangle$, then observes that 
$M$ corresponds to a complex subvariety $Z$
of $\C^N$, smooth outside of $\{0\}$ and fixed by $A$.
The operator $A$ admits a Jordan-Chevalley decomposition $A:= S U$,  with $S$ diagonal and 
$U$ unipotent and one can show that $S$ preserves  $Z$.
Then $M_1 := (Z\setminus \{0\})/\langle S\rangle$ is a deformation of $M$ (as $S$ is
contained in
a  $\mathrm{GL}(\C^n)$-orbit of $A$) 
and is Vaisman as {contained in the Hopf manifold $H_S:= (\C^n\setminus \{ 0\})/\langle S\rangle$.}

\smallskip

The above result shows that all known topological obstructions to the existence of a Vaisman metric on a compact complex manifold (see {\em e.g.} \cite{drag}) apply to LCK manifolds with potential. It allows, in particular, to determine the fundamental group of compact LCK manifolds with potential. Indeed, one first deforms the structure to a Vaisman one, then deforms this one to a quasi-regular one (see above) which fibers in elliptic curves over a K\"ahler basis $X$. At this point, one considers the homotopy sequence of the fibering:
$$\pi_2(X) \stackrel \delta
\arrow \pi_1(T^2) \arrow \pi_1(M) \arrow \pi_1(X) \arrow 0$$  
and observes that $\rk(\Im(\delta))\leq 1$ in $\pi_1(T^2)$, as the Chern classes of the $S^1\times S^1$-fibration
 are: one trivial (as $M$ fibers on $S^1$),
the other one non-trivial,
as $M$ is non-K\"ahler, and the total space of
an elliptic fibration with trivial Chern classes is K\"ahler. Hence:

\begin{corollary}\cite{_OV_Top_Potential_}
 The fundamental group of a
compact LCK manifold $M$ with an automorphic potential
admits an exact sequence
\[
0 \arrow G \arrow \pi_1(M) \arrow \pi_1(X) \arrow 0
\]
where $\pi_1(X)$ is the fundamental group of a K\"ahler  orbifold,
and $G$ is a quotient of $\Z^2$ by a subgroup of rank 1.
\end{corollary}

\begin{remark}
 In fact, in \cite{_OV_Top_Potential_} we only proved that the rank of the subgroup must be $\leq 1$, but the recent Theorem \ref{vuli} above (\cite{vu2}) shows that $\rk(M)=0$ would imply $M$ is K\"ahler (see also Remark \ref{rem1}).
\end{remark}

\begin{corollary}\cite{_OV_Top_Potential_}
 A non-Abelian free group cannot be the
fundamental group of a compact LCK manifold
with potential.
\end{corollary}

This corollary, as well as other topological restrictions, was first obtained by Kokarev and Kotschick using harmonic forms and a LCK version of Siu-Beauville result:

\begin{theorem}\cite{_Kokarev_Kotschick:Fibrations_}
 Let $M$ be a closed complex manifold admitting a LCK structure with potential (pluricanonical K\"ahler-Weyl). Then
the following statements are equivalent:

$i)$ $M$ admits a surjective holomorphic map with connected fibers to a closed Riemann surface
of genus $\geq 2$;

$ii)$ $\pi_1(M)$ admits a surjective homomorphism to the fundamental group
of a closed Riemann surface of genus $\geq 2$;

$iii)$ $\pi_1(M)$  admits a surjective homomorphism to a non-Abelian free
group.
\end{theorem}

The above can be generalized to:

\begin{theorem}\cite{_Kokarev_Kotschick:Fibrations_}
 Let $M$ be a closed complex manifold admitting a LCK structure with potential, and
$N$ a closed Riemannian manifold of constant negative curvature. If $\f: \pi_1(M)\rightarrow \pi_1(N)$ is
a representation with non-cyclic image, then there exists a compact Riemann surface $S$ and a
holomorphic map $h: M \rightarrow S$ with connected fibers such that $\f$ factors through $h_*$.

In particular, if $N$ is a closed real hyperbolic
manifold, $\dim N\geq 4$, then any  map $f: M\rightarrow N$ has degree zero.
\end{theorem}

Other topological obstructions to the existence of a LCK structure with potential were obtained by Kokarev in \cite{kok} using harmonic maps techniques. For example, one of his results is:

\begin{theorem}\cite{kok}
 Let $M$ be a compact LCK manifold of the same homotopy
type as a locally Hermitian symmetric space of non-compact type whose universal cover does
not contain the hyperbolic plane as a factor. If $M$ admits a LCK metric with potential, then it
admits a global K\"ahler metric.
\end{theorem}

On the other hand, on compact Vaisman manifolds the cohomology of $L$ (which is the Morse-Novikov cohomology of the operator $d-\theta\wedge$) is simple: $H^*(M,L_\theta)=0$ follows easily from the Structure theorem \ref{ov_str} (here the subscript $\theta$ makes precise the structure of local system of $L$). 

\begin{theorem}\cite{_ov:MN_}\label{mnvai}
 Let $(M J)$ be a compact complex manifold, of complex dimension at least $3$, endowed with a Vaisman structure with $2$-form $\omega$  and Lee form $\theta$. Let $\omega_1$ be another LCK-form (not necessarily Vaisman) on $(M J)$, and let $\theta_1$ be its Lee form. Then $\theta_1$ is cohomologous with
the Lee form of a Vaisman metric, and $[\omega_1]=0\in H^2(M,L_{\theta_1})$.
\end{theorem}

By contrast, on an Inoue surface, which does not admit any Vaisman metric, there exists a LCK metric, compatible with the solvmanifold structure, with non-vanishing Morse-Novikov class of the LCK two-form, \cite{ba}.

\medskip

We end this section with a result which determines all compact nilmanifolds admitting an invariant LCK structure (generalizing a result of L. Ugarte in dimension 4):

\begin{theorem}\cite{sawai}
 Let $(M, J )$ be a non-toral compact nilmanifold with a left-invariant complex
structure. If $(M, J )$ has a locally conformally K\"ahler structure, then $(M, J )$ is biholomorphic
to a quotient of $(H(n) \times \R, J_0)$, where $H(n)$ is the generalized Heisenberg group and $J_0$ is the natural complex structure on the product.
\end{theorem}
 
The author mentions that he does not know if the biholomorphism he finds passes to the quotient; in other words, he does not know if the compact LCK nilmanifold is isomorphic or biholomorphic with $H(n)\times S^1$. On the other hand, one sees that, in particular, left invariant LCK structures on compact nilmanifolds are of Vaisman type. We tend to believe that the result is true in more general setting, namely without the assumption of left (or right) invariance. It is tempting to state:

\begin{conjecture}
 Every LCK compact nilmanifod is, up to covering, the product of the generalized Heisenberg group with $S^1$.
\end{conjecture}

\section{Transformation groups of LCK manifolds}

The study of this topic went in two directions. The first one is characterizing the various groups appearing in LCK geometry (conformalities, isometries, affinities with respect to the Levi-Civita or the Weyl connection, holomorphicities) and determination of their interrelations. The second one is characterizing different subclasses of LCK manifolds by the existence of a particular subgroup of one of these groups.

In the first direction,  we mention the above Theorem \ref{conf_hol} and the following {\em local} result:

\begin{theorem}\cite{mo}
 On any LCK manifold, $\mathfrak{aff}(M,\nabla)=\mathfrak{aut}(M)$, provided that
   $\mathrm{Hol}_0(D)$ is irreducible and $\mathrm{Hol}_0(D)$ is not contained in $\mathrm{Sp}(n)$.
\end{theorem}
 
For the proof, one first shows that $\mathfrak{aff}(M,\nabla)\subseteq \mathfrak{h}(M,J)$ (thus generalizing the analogue result for K\"ahler manifolds). 
Indeed, let $f\in  \mathrm{Aff}(M,D)$. We show that it is $\pm$ - holomorphic.  

Define $J'_x:=(d_xf)^{-1}\circ
J_{f(x)}\circ (d_xf)$. Then $J'$ is $D$-parallel.
To show that $J'=\pm J$, we 
decompose $JJ'=S\,(\text{symm.})+A\,(\text{antisymm.})$.
Then $S$ is $\nabla$-parallel and hence it has constant eigenvalues; thus the  
corresponding eigenbundles are $D$-parallel.

By $\mathrm{Hol}_0(D)$ irreducible, $S=k\id$, $k\in\mathbb{R}$.
Similarly,  $A^2=p\id$, $p\in\mathbb{R}$.  

Now, if $A\neq 0$, then $A(X) \ne 0$ for some $X\in TM$, so
$0>-g(AX,AX)=g(A^2X,X)=pg(X,X),$ whence $p<0$.
Then $K:=A/\sqrt{-p}$ is $D$-parallel, $K^2=-\id$, 
$KJ=-JK$, so $(J,K)$ defines a $D$-parallel quaternionic structure 
structure on $M$, contradiction. 

Hence $A=0$, $JJ'=k\id$, so $J'=-kJ$.
But ${J'} ^2=-\id$, thus $k=\pm 1$ and so
$J'=\pm J$.

With similar arguments one proves that $\mathfrak{aff}(M,D)\subseteq \mathfrak{conf}(M,[g])$.

\smallskip

In the second direction, we first recall the following characterization of Vaisman manifolds:

\begin{theorem}\cite{_Kamishima_Ornea_}\label{kamo}
A compact LCK manifold admits a LCK
metric with parallel Lee form if its Lie group of holomorphic
conformalities has  a complex one-dimensional Lie subgroup,
acting non-isometrically on its K\"ahler covering.
\end{theorem}

We note that the above criterion assures the existence of a Vaisman metric {\em in the conformal class} of the given LCK one. We recently extended this result to obtain the existence of a LCK metric with automorphic potential, not necessarily conformal to the starting LCK one:

\begin{theorem}\cite{_OV_Transf_}
Let $M$ be a compact LCK manifold, equipped
with a holomorphic $S^1$-action. Suppose that
the holonomy of the weight bundle $L$ restricted to a general
orbit of this $S^1$-action is non-trivial. Then $\tilde M$ admits a global
automorphic potential. 
\end{theorem}

For the proof, a first step is to show that one can assume from the beginning  that $\omega$, and hence, as $J$ remains unchanged,  $g$, is $S^1$-invariant ({\em i.e.} the action is isometric). Note that a similar argument was used in the proof of Theorem \ref{mnvai}.

Indeed, we average $\theta$ on $S^1$ and obtain  $\theta'=\theta+df$ which is   $S^1$-invariant. The cohomology class is conserved:  $[\theta]=[\theta']$.
Now we let 
$\omega'=e^{-f}\omega$: it is LCK, conformal
to $\omega$ and with Lee form $\theta'$. 

Hence, we may assume from the beginning that
$\theta$ (corresponding to $\omega$) is $S^1$-invariant.

We now average $\omega$ over $S^1$, taking into account that:
\begin{equation}
   d(a^*\omega) = a^*\omega\wedge a^*\theta
   = a^*\omega\wedge \theta, \quad\quad a\in S^1.
\end{equation}
We thus find an 
 $\omega'$ which is $S^1$-invariant, with 
\[
d\omega' = \theta\wedge\omega'.
\]

As the monodromy of $L$ along
an orbit $S$ of the $S^1$-action is precisely $\int_S \theta$, it is 
not changed
by this averaging procedure.

This means that it is enough to make the proof 
assuming $\omega$ is $S^1$-invariant.

On the other hand, the lift of $S^1$
to $\tilde M$ acts on $\tilde \omega$ by homotheties, and the corresponding
conformal constant is equal to the
monodromy of $L$ along the orbits of
$S^1$. Thus, the image of the restriction of the character $\chi$ to the lifted subgroup cannot be compact in $\R^{>0}$ unless it is trivial, hence the $S^1$ action lifts to an $\R$ action.

In conclusion, we may assume from the beginning that
$S^1$ is lifted to an $\R$ acting on
$\tilde M$ by non-trivial homotheties.

Now, let $A$ be the vector field on
$\tilde M$ generated by the
$\R$-action. $A$ is holomorphic and homothetic ($\Ll_A\Omega  =\la\Omega$). 

Let $A^c=JA.$ Then:
$$dd^c|A|^2= \la^2\Omega  +\Ll_{A^c}^2\Omega$$
Read in Bott-Chern cohomology, this implies:
$$\Ll_{A^c}^2[\Omega]=-\la^2[\Omega]\in H^2_{BC}(M,L).$$
 Hence
$V:=\mathrm{span}\{[\Omega], \Ll_{A^c}[\Omega]\}\subset H^2_{BC}(M,L)$ is
$2$-dimensional. 

As $\Ll_{A^c}$ acts on $V$ with two
1-dimensional eigenspaces, corresponding
to $\1\lambda$ and $-\1\lambda$, it is essentially a rotation with $\lambda\pi/2$, and hence the flow of $A^c$ 
satisfies:
$$e^{t A^c}[\Omega]=[\Omega], \, \mathrm{for}\,\, t=
2n\pi\lambda^{-1},\, n\in\Z.$$
But also
$$
 \int_0^{2\pi\lambda^{-1}} e^{t A^c}[\Omega]dt=0.
$$
and hence 
$$\Omega_1:= \int_0^{2\pi\lambda^{-1}} e^{t A^c}\Omega dt$$
is a {K\"ahler  form}, whose Bott-Chern class vanishes, $[\Omega_1]=0\in H^2_{BC}(M,L)$, thus defining a LCK metric with automorphic potential by Lemma \ref{potbc}.



\end{document}